\documentclass[10pt,twoside]{article}
\usepackage{amsmath}
\usepackage{amsthm}
\usepackage{amssymb}
\usepackage{Latex-document}

\def\bysame{\leavevmode\hbox to3em{\hrulefill}\thinspace}

\renewcommand{\labelenumi}{\theenumi}
\renewcommand{\theenumi}{\rm (\roman{enumi})}

\newtheorem{thm}{Theorem}

\newtheorem{prop}[thm]{Proposition}
\newtheorem{conj}[thm]{Conjecture}
\newtheorem*{question}{Question}
\theoremstyle{remark}
\newtheorem*{pro}{\bf Proof}
\newtheorem*{rem}{\bf Remark}
\newtheorem*{exa}{\bf Example}
\theoremstyle{definition}
\newtheorem{defn}[thm]{Definition}

\newfont{\cyrr}{wncyr10}
\def\Sh{\mbox{\cyrr Sh}}

\def\Z{\mathbf{Z}}
\def\Q{\mathbf{Q}}

\def\I{\mathbf{I}}
\def\L{\mathbf{L}}

\def\C{\mathbf{C}}

\def\bLambda{\mathbf{\Lambda}}
\def\bGamma{\mathbf{\Gamma}}

\def\cA{J}
\def\N{{\mathcal{N}}}
\def\O{\mathcal{O}}
\def\cL{\mathcal{L}}
\def\OK{\O_K}

\def\Gal{\mathrm{Gal}}
\def\Sym{\mathrm{Sym}}
\def\Hom{\mathrm{Hom}}

\def\cycl{\mathrm{cycl}}
\def\anti{\mathrm{anti}}
\def\rk#1{\mathrm{rank}\,#1}
\def\Sel{\mathrm{Sel}}
\def\cp{\mathrm{char}}
\def\tors{\mathrm{tors}}
\def\disc{\mathrm{disc}}

\def\Zp{\Z_p}
\def\Qp{\Q_p}
\def\Qb{\bar{\Q}}
\def\Kb{\bar{K}}

\def\Kcycl{K_\infty^{\cycl}}
\def\Kanti{K_\infty^{\anti}}
\def\Ktwo{\mathbf{K}_\infty}
\def\opp{(\tau)}
\def\map#1{\;\xrightarrow{#1}\;}
\def\isom{\map{\sim}}
\def\medoplus{\mathop{\hbox{\Large$\oplus$}}}
\def\Selp#1{\mathcal{S}_p(E_{/#1})}
\def\trivchar{{\chi_{\mathrm{triv}}}}
\def\H{\mathcal{H}}
\def\U{\mathcal{U}}
\def\h{\mathrm{h}}
\def\Ganti{\Gamma_{\hskip -2pt\anti}}
\def\Gcycl{\Gamma_{\hskip -2pt\cycl}}

\title{\bf  Elliptic Curves and Class Field Theory\vskip 6mm}

\author{Barry Mazur\thanks{Department of Mathematics,
Harvard University, Cambridge, MA 02138, USA. E-mail:
mazur@math.harvard.edu } \quad {\bf Karl Rubin}\thanks{Department
of Mathematics, Stanford University, Stanford, CA 94305, USA.
E-mail:\hskip 14mm rubin@math.stanford.edu } \vspace*{-0.5cm}}
\date{\vspace{-8mm}}

\begin{document}
\maketitle

\thispagestyle{first} \setcounter{page}{185}

\begin{abstract}\vskip 3mm

Suppose $E$ is an elliptic curve defined over $\Q$.
At the 1983 ICM the first author formulated some
conjectures that propose a close relationship
between the explicit class field theory construction of certain
abelian extensions of imaginary quadratic fields and an explicit
construction that (conjecturally) produces almost all of the
rational points on $E$ over those fields.

Those conjectures are to a large extent settled by recent work of
Vatsal and of Cornut, building on work of Kolyvagin and others.
In this paper we describe a collection of interrelated
conjectures still open regarding the variation of
Mordell-Weil groups of $E$ over abelian extensions of
imaginary quadratic fields, and suggest a possible
algebraic framework to organize them.

\vskip 4.5mm


\end{abstract}

\vskip 12mm

\section{Introduction} \label{intro}\setzero
\vskip-5mm \hspace{5mm }

Eighty years have passed since Mordell proved that the (Mordell-Weil) group of
rational points on an elliptic curve $E$ is finitely generated, yet
so limited is our knowledge that
we still have no algorithm guaranteed to compute the rank of this group.
If we want to ask even more ambitious questions about how the rank of
the Mordell-Weil group $E(F)$ varies as $F$ varies, it
makes sense to restrict attention only to those fields for which we
have an explicit construction, such as finite abelian
extensions of a given imaginary quadratic field $K$.  Taking our lead from
the profound discovery of Iwasawa that the variational properties
of certain arithmetic invariants are well-behaved if one restricts to
subfields of $\Zp^d$-extensions of number fields, we will focus on the
following Mordell-Weil variation problem:

\begin{quote}
\sl
Fixing an elliptic curve $E$ defined over $\Q$, an imaginary quadratic field $K$,
and a prime number $p$, study the variation of the Mordell-Weil group of
$E$ over finite subfields of the (unique) $\Zp^2$-extension of $K$ in $\Kb$.
\end{quote}

This problem was the subject of some conjectures formulated by the first
author at the 1983 ICM \cite{mazuricm}, conjectures which have recently been largely
settled by work of Vatsal \cite{vatsal} and Cornut \cite{cornut} building on
work of Kolyvagin and others.

\begin{exa}
Let $E$ be the elliptic curve $y^2 +y = x^3 - x$, $p = 5$, and let $K = \Q(\sqrt{-7})$.
If $F$ is a finite extension of $K$, contained in the $\Z_5^2$ extension of $K$,
then $\rk{E(F)} = [F \cap \Kanti : K]$ where $\Kanti$ is the anticyclotomic
$\Z_5$-extension of $K$ (see \S2 for the definition).
One only has an answer like this in the very simplest cases.

Now with the same $E$ and $p$, take $K = \Q(\sqrt{-26})$.  A guess here would be
that $\rk{E(F)} = [F \cap \Kanti : K] + 2$, but this seems to be
beyond present technology.
\end{exa}

The object of this article is to sketch a package of still-outstanding
conjectures in hopes that it offers an even more
precise picture of this piece of arithmetic.  These conjectures are in some
cases due to, and in other cases build on ideas of, Bertolini \& Darmon,
Greenberg, Gross \& Zagier, Haran, Hida, Iwasawa, Kolyvagin, Nekov\'a\u{r},
Perrin-Riou, and the authors, among others.

In sections 3 through 5 we describe the three parts of our picture:
the {\em arithmetic theory} (the study of the
Selmer modules over Iwasawa rings that contain the information we seek),
the {\em analytic theory} (the construction and
study of the relevant $L$-functions, both classical and $p$-adic),
and the {\em universal norm theory} which arises from
purely arithmetic considerations, but provides analytic
invariants.

In the final section we suggest the beginnings of a new algebraic structure
to organize these conjectures.  This structure should not be viewed
as a conjecture, but rather as a mnemonic to collect our
conjectures and perhaps predict new ones.

More details and proofs will appear in a forthcoming paper.

\section{Running hypotheses and notation} \label{notation}\setzero
\vskip-5mm \hspace{5mm }

Fix a triple $(E, K, p)$  where $E$ is an elliptic curve of conductor
$N$ over $\Q$, $K$ is an imaginary quadratic field of discriminant $D < -4$, and
$p$ is a prime number.  To keep our discussion focused and as succinct as
possible, we make the following hypotheses and conventions.

Assume that $p$ is odd, that $N$, $p$ and $D$
are pairwise relatively prime, and that if $E$ has complex multiplication,
then $K$ is {\em not} its field of complex multiplication. Let
$\OK \subset K$ denote the ring of integers of $K$.
Assume further that there exists an ideal $\N \subset \OK$
such that $\OK/\N$ is cyclic of order $N$ (this is sometimes called the
{\bf Heegner Hypothesis}), and that $p$ is a prime of ordinary reduction
for $E$.
For simplicity we will assume throughout this article that the $p$-primary
subgroups of the Shafarevich-Tate groups of $E$ over the number fields we
consider are all finite.

\begin{prop}
\label{parity}
Under the assumptions above, $\rk{E(K)}$ is odd.
\end{prop}

\begin{pro}
This follows from the Parity Conjecture recently proved by
Nekov\'a\u{r} \cite{nekovar}.
\end{pro}

Let $\Ktwo$ denote the (unique) $\Zp^2$-extension of $K$ and
$\bGamma := \Gal(\Ktwo/K)$, so $\bGamma \cong \Zp^2$.  We define the Iwasawa ring
$$
\bLambda := \Zp[[\bGamma]] \otimes_{\Zp} \Qp.
$$
(To simplify notation and to avoid some complications, we will often
work with $\Qp$-vector spaces instead of natural $\Zp$-modules;
in particular we have tensored the usual Iwasawa ring with $\Qp$.)
For every (finite or infinite) extension F of K in $\Ktwo$ we also define
$$
\Lambda_F := \Zp[[\Gal(F/K)]] \otimes_{\Zp} \Qp,
    \qquad \I_F := \ker\{\bLambda \twoheadrightarrow \Lambda_F\}.
$$
Then $\I_K$ is the augmentation ideal of $\bLambda$,
and if $[F:K]$ is finite then $\Lambda_F$ is just the group ring $\Qp[\Gal(F/K)]$.
If $\Gal(F/K)$ is $\Zp$ or $\Zp^2$, and $M$ is a finitely generated
torsion $\Lambda_F$-module, then $\cp_{\Lambda_F}(M)$ will denote
the characteristic ideal of $M$.  In
particular $\cp_{\Lambda_F}(M)$ is a principal ideal of $\Lambda_F$.

There is a $\Q_p$-projective line of $\Zp$-extensions of $K$, all
contained in $\Ktwo$.  Among these are two distinguished $\Zp$-extensions:
\begin{itemize}
\item
the {\bf cyclotomic} $\Zp$-extension $\Kcycl$,
the compositum of $K$ with the unique (cyclotomic) $\Zp$-extension of $\Q$
(write $\Gcycl = \Gal(\Kcycl/K)$, $\Lambda_{\cycl} = \Lambda_{\Kcycl}$),
\item
the {\bf anticyclotomic} $\Zp$-extension $\Kanti$, the
unique $\Zp$-extension of $K$ that is Galois over $\Q$ with non-abelian,
and in fact dihedral, Galois group (write $\Ganti = \Gal(\Kanti/K)$,
$\Lambda_{\anti} = \Lambda_{\Kanti}$).
\end{itemize}
Then $\bGamma = \Gcycl \oplus \Ganti$ and
$\bLambda = \Lambda_{\cycl} \otimes_{\Zp} \Lambda_{\anti}$.

Complex conjugation $\tau : K \to K$ acts on $\bGamma$, acting as $+1$ on
$\Gcycl$ and $-1$ on $\Ganti$.  This induces nontrivial
involutions of $\bLambda$ and $\Lambda_{\anti}$, which we also
denote by $\tau$.
If $M$ is a module over $\bLambda$ (or similarly over $\Lambda_{\anti}$),
let $M^{\opp}$ denote the module whose underlying abelian group is $M$ but
where the new action of $\gamma\in \bGamma$ on $m \in M^{\opp}$ is given by
the old action of  $\gamma^{\tau}$ on $m$.

Our $\bLambda$-modules will usually come with a natural action
of $\Gal(\Ktwo/\Q)$.  These actions are
continuous and $\Zp$-linear, and satisfy the formula
$\tilde\tau(\gamma\cdot m) =\gamma^{\tau}\cdot\tilde\tau(m)$
for every lift $\tilde\tau$ of $\tau$ to $\Gal(\Ktwo/\Q)$.
Thus the action of any lift $\tilde\tau$ induces an isomorphism $M \isom M^{\opp}$.
We will refer to such $\bLambda$ or $\Lambda_{\anti}$-modules as
{\bf semi-linear $\tau$-modules}.  If $M$ is a semi-linear $\tau$-module
and is free of rank one over $\Lambda_{\anti}$, we define the {\bf sign}
of $M$ to be the sign $\pm 1$ of the action of $\tau$ on the one-dimensional
$\Qp$-vector space $M\otimes_{\Lambda_{\anti}}{\Lambda_K}$.
Such an $M$ is completely determined (up to
isomorphism preserving its structure) by its sign.

\begin{defn}
\label{herm}\it If $M$ and $A$ are semi-linear $\tau$-modules,
then a ($\bLambda$-bilinear) $A$-valued {\bf $\tau$-Her\-mi\-tian
pairing} $\pi$ is a $\bLambda$-module homomorphism $ \pi : M
\otimes_{\bLambda} M^{\opp} \to A $ such that for every lift
$\tilde\tau$ of $\tau$ to $\Gal(\Ktwo/\Q)$
$$
\pi(m \otimes n) = \pi(n \otimes m)^{\tilde\tau} = \pi(\tilde\tau n \otimes \tilde\tau m).
$$\rm
\end{defn}

\section{Universal norms} \label{universalnorms}\setzero
\vskip-5mm \hspace{5mm }

\begin{defn}\it
If $K \subset F \subset \Ktwo$, the {\bf universal norm module}
$U(F)$ is the projective limit
$$
U(F) := \Qp \otimes \varprojlim_{K \subset L \subset F}(E(L)\otimes \Zp)
$$
(projective limit with respect to traces, over {\em finite}
extensions $L$ of $K$ in $F$) with its natural
$\Lambda_F$-structure.  If $F$ is a finite extension of $K$, then
$U(F)$ is simply $E(F) \otimes \Qp$.\rm
\end{defn}

If $F$ is a $\Zp$-extension of $K$, then $U(F)$ is a free $\Lambda_F$-module
of finite rank, and is zero if and only if the Mordell-Weil ranks of $E$
over subfields of $F$ are bounded
(cf.\ \cite{mazuricm} \S18 or \cite{perrinrioubsmf} \S2.2).
The first author conjectured some time ago \cite{mazuricm} that
for $\Zp$-extensions $F/K$, and under our running hypotheses,
$U(F) = 0$ if $F \ne \Kanti$ and $U(\Kanti)$ is free of rank
one over $\Lambda_{\anti}$.  The following theorem follows from
recent work of Kato \cite{kato}
for $\Kcycl$ and Vatsal \cite{vatsal} and Cornut \cite{cornut}
for $\Kanti$.

\begin{thm}
\label{growthnumber}
$U(\Kcycl) = 0$ and $U(\Kanti)$ is free of rank one over $\Lambda_{\anti}$.
\end{thm}

For the rest of this paper we will write $\U$ for the anticyclotomic
universal norm module $U(\Kanti)$.  Complex conjugation gives $\U$ a natural
semi-linear $\tau$-module structure.  Since $\U$ is free of rank one over
$\Lambda_{\anti}$, we conclude that $\U$ is completely determined
(up to isomorphism preserving its $\tau$-structure) by its sign.

Let $r^\pm$ be the rank of the $\pm1$ eigenspace of $\tau$ acting on $E(K)$,
so $\rk{E(\Q)} = r^+$ and $\rk{E(K)} = r^+ + r^-$.
By Proposition \ref{parity}, $\rk{E(K)}$ is odd so $r^+ \ne r^-$.

\begin{conj}[Sign Conjecture]
\label{signconj}
The sign of the semi-linear $\tau$-module $\U$ is
$+1$ if $r^+ > r^-$, and is $-1$ if $r^- > r^+$.
\end{conj}

\begin{rem}
Equivalently, the Sign Conjecture asserts that the sign of $\U$ is
$+1$ if twice $\rk{E(\Q)}$ is greater than $\rk{E(K)}$, and $-1$ otherwise.

As we discuss below in \S4, the Sign Conjecture is related to
the nondegeneracy of the $p$-adic height pairing (see the remark
after Conjecture \ref{maxnondeg}).
\end{rem}

The $\Lambda_{\anti}$-module $\U$ comes with a canonical Hermitian
structure.  That is, the canonical (cyclotomic) $p$-adic height pairing
(see \cite{mazurtate} and \cite{perrinrioubsmf} \S2.3)
$$
\h : \U \otimes_{\Lambda_{\anti}}\U^{\opp} \longrightarrow
    \Gcycl\otimes_{\Zp}\Lambda_{\anti}
$$
is a $\tau$-Hermitian pairing in the sense of Definition \ref{herm}.

\begin{conj}[Height Conjecture]
\label{heightconj}
The homomorphism $\h$ is an isomorphism of free $\Lambda_{\anti}$-modules of rank one
$$
\h : \U \otimes_{\Lambda_{\anti}}\U^{\opp} \isom \Gcycl\otimes_{\Zp}\Lambda_{\anti}.
$$
\end{conj}

The $\Lambda_{\anti}$-module $\U$ has an important submodule, the {\bf Heegner submodule} $\H \subset  \U$. Fix a
modular parameterization $X_0(N) \to E$. The Heegner submodule $\H$ is the cyclic $\Lambda_{\anti}$-module
generated by a trace-compatible sequence $c =\{c_L\}$ of Heegner points $c_L \in E(L)\otimes \Zp$ for finite
extensions $L$ of $K$ in $\Kanti$. See for example \cite{mazuricm} \S19 or \cite{perrinrioubsmf} \S3.  Call such a
$c \in \H$ a {\bf Heegner generator}.  The Heegner generators of $\H$ are well-defined up to multiplication by an
element of $\pm \Gamma \subset (\Lambda_{\anti})^\times$. The $\Lambda_{\anti}$-submodule $\H\subset \U$ is stable
under the semi-linear $\tau$-structure of $\U$, so the action of $\tau$ gives an isomorphism $\U/\H \isom
(\U/\H)^{\opp} \cong \U^{\opp}/\H^{\opp}$.

Let $c^{\opp}$ denote the element $c$ viewed in the
$\Lambda_{\anti}$-module $\H^{\opp}$.  Since
$$
(\pm \gamma c)\otimes_{\Lambda_{\anti}} (\pm \gamma c)^{\opp}
    = c \otimes_{\Lambda_{\anti}} c^{\opp}
$$
for every $\pm \gamma \in \pm \Gamma$, the element
$c \otimes c^{\opp} \in \H \otimes_{\Lambda_{\anti}}\H^{\opp}$
is independent of the choice of Heegner generator, and is therefore a totally
canonical generator of the free, rank one $\Lambda_{\anti}$-module
$\H \otimes_{\Lambda_{\anti}}\H^{\opp}$.

\begin{defn}\it
The {\bf Heegner $L$-function} (for the triple
$(E,K,p)$ satisfying our running hypotheses) is the element
$$
\cL \;:=\; \h(c \otimes c^{\opp})
    \in \Gcycl\otimes_{\Zp}\Lambda_{\anti}.
$$\rm
\end{defn}

\begin{conj}
\label{heegnerlfunction}
\quad
$\Gcycl \otimes \cp(\U/\H)^2 = \Lambda_{\anti}\cL$ \quad
inside $\Gcycl\otimes \Lambda_{\anti}$.
\end{conj}

One sees easily that
$\Gcycl \otimes \cp(\U/\H)^2 \supset \Lambda_{\anti}\cL$,
and that Conjecture \ref{heegnerlfunction}
is equivalent to the Height Conjecture (Conjecture \ref{heightconj}).

\section{ The analytic theory} \label{analytic}\setzero
\vskip-5mm \hspace{5mm }

The (``two-variable") $p$-adic $L$-function for $E$ over $K$ is an element
$\L \in \bLambda$
constructed by Haran \cite{haran} and by a different,
more general, method by Hida \cite{hida} (see also the papers of
Perrin-Riou \cite{perrinriouinv, perrinrioujlms}).
The $L$-function $\L$ is characterized by the fact that it interpolates
special values of the classical Hasse-Weil $L$-function of twists of
$E$ over $K$.  More precisely, embedding $\Qb$ both in $\C$ and
$\Qb_p$,  if
$
\chi:  \bGamma \to {\bar{\Z}}^\times \subset \bar{\Z}_p^\times
$
is a character of finite order then
\begin{equation}
\label{const}
\chi(\L) = c(\chi)\frac{L_{\rm classical}(E_{/K},\chi,1)}{8 \pi^2 \|f_E\|^2}
\end{equation}
where $L_{\rm classical}(E_{/K},\chi,s)$ is the Hasse-Weil $L$-function of
the twist of $E_{/K}$ by $\chi$,
$c(\chi)$ is an explicit algebraic number (cf.\ \cite{perrinriouinv} Th\'eor\`eme 1.1),
$f_E$ is the modular form on $\Gamma_0(N)$
corresponding to $E$, and $\|f_E\|$ is its Petersson norm.

Projecting $\L\in\bLambda$ to the cyclotomic or the
anticyclotomic line
via the natural projections $\bLambda \to \Lambda_{\cycl}$ and
$\bLambda \to \Lambda_{\anti}$, we get ``one-variable" $p$-adic $L$-functions
$$
\L \mapsto L_{\cycl} \in\Lambda_{\cycl}
\qquad\text{and}\qquad
\L \mapsto L_{\anti} \in \Lambda_{\anti}.
$$
It follows from the functional equation satisfied by $\L$
(\cite{perrinriouinv} Th\'eor\`eme 1.1)
and the Heegner Hypothesis that $L_{\anti} = 0$. In other words, viewing
$\bLambda = \Lambda_{\anti}[[\Gcycl]]$
as the completed group ring of $\Gcycl$ with coefficients
in $\Lambda_{\anti}$, we have that the
``constant term'' of $\L \in \Lambda_{\anti}[[\Gcycl]]$
vanishes.  We now consider its ``linear term."

There is a canonical isomorphism of (free, rank one) $\Lambda_{\anti}$-modules
$$
\Gcycl\otimes_{\Zp}\Lambda_{\anti} \cong \I_{\Kanti}/\I_{\Kanti}^2
$$
which sends
$\gamma\otimes 1 \in \Gcycl\otimes_{\Zp} \Lambda_{\anti}$
to $\gamma -1 \in \I_{\Kanti}/\I_{\Kanti}^2$.

\begin{conj}[$\Lambda$-adic Gross-Zagier Conjecture]
\label{lambdagz}
Let $L'$ denote the image of $\L$ under the map
$\I_{\Kanti}/\I_{\Kanti}^2 \isom \Gcycl\otimes_{\Zp}\Lambda_{\anti}$.
Then
$$
L' = d^{-1}\cL
$$
where $d$ is the degree of the modular parametrization $X_0(N) \to E$.
\end{conj}

\begin{rem}
Perrin-Riou \cite{perrinriouinv} proved that if $p$ splits in $K$ and
the discriminant $D$ of $K$ is odd, then
$L'$ and $d^{-1}\cL$ have the same image under the projection
$
\Lambda_{\anti} \to \Lambda_K = \Qp.
$
\end{rem}

Let  $\I := \I_K$, the augmentation ideal of $\bLambda$.
For every integer $r \ge 0$ we have $\I^r/\I^{r+1}\cong \Sym_{\Zp}^r(\bGamma)\otimes\Qp.$
Using the direct sum decomposition
$\bGamma = \Gcycl \oplus \Ganti$
we get a canonical direct sum decomposition
\begin{equation}
\label{decomp}
\Sym_{\Zp}^r(\bGamma) = \medoplus_{j=0}^r \Gamma^{r-j,j}
\quad \text{where
$\Gamma^{i,j} := (\Gcycl)^{\otimes i} \otimes_{\Zp} (\Ganti)^{\otimes j}.$}
\end{equation}

Consider the canonical (two-variable) $p$-adic height pairing
\begin{equation}
\label{2varheight}
\langle ~,\;\rangle \;:\; E(K) \times E(K) \longrightarrow \bGamma\otimes\Qp.
\end{equation}
Set $r = \rk{E(K)}$, which is odd by Proposition \ref{parity}.
Define the {\bf two-variable $p$-adic regulator}
$R_p(E,K)$ to be the discriminant of this pairing:
$$
R_p(E,K) := t^{-2} \det \langle P_i,P_j\rangle \in
\Sym_{\Zp}^r(\bGamma)\otimes\Qp \cong \I^r/\I^{r+1},
$$
where $\{P_1,\ldots,P_r\}$ generates a subgroup of $E(K)$ of finite index $t$.
 For each integer $j = 0, \dots , r$ let $R_p(E,K)^{r-j,j}$ be the projection
of $R_p(E,K)$ into $\Gamma^{r-j,j}\otimes\Qp$ under \eqref{decomp}, so that
$$
R_p(E,K) = \medoplus_{j=0}^r \ R_p(E,K)^{r-j,j}.
$$

Recall that $r^{\pm}$ is the rank of the $\pm1$-eigenspace
$E(K)^{\pm}$ of $\tau$ acting on $E(K)$.

\begin{prop}
\label{deg}
 $R_p(E,K)^{r-j,j} = 0$ unless $j$ is even and $j \le 2 \min(r^+, r^-)$.
\end{prop}

\begin{pro}
This follows from the fact that the height pairing \eqref{2varheight}
is $\tau$-Hermitian, so
$\langle \tau x, \tau y \rangle = \langle x, y \rangle^{\tau}$,
and therefore the induced height pairings
$$
E(K)^\pm \times E(K)^\pm \to \Ganti\otimes\Qp, \quad
E(K)^+ \times E(K)^- \to \Gcycl\otimes\Qp
$$
vanish.
\end{pro}

\begin{conj}[Maximal nondegeneracy of the height pairing]
\label{maxnondeg}
If $j$ is even and $0 \le j \le 2 \min(r^+, r^-)$ then
$R_p(E,K)^{r-j,j} \ne 0$.
\end{conj}

\begin{rem}
Conjecture \ref{maxnondeg}, or more specifically the nonvanishing of
$R_p(E,K)^{r-j,j}$ when $j = 2 \min(r^+, r^-)$, implies the Sign Conjecture
(Conjecture \ref{signconj}).  This is proved in the same way as
Proposition \ref{deg}, using the additional fact that the anticyclotomic universal norms in
$E(K) \otimes \Zp$ are in the kernel of the anticyclotomic $p$-adic height pairing
$(E(K) \otimes \Zp) \times (E(K) \otimes \Zp) \to \Ganti\otimes\Qp$.
\end{rem}

\section{The arithmetic theory} \label{arithmetic}\setzero
\vskip-5mm \hspace{5mm }

For every algebraic extension $F$ of $K$, let $\Sel_p(E_{/F})$ denote
the $p$-power Selmer group of $E$ over $F$, the subgroup of $H^1(G_F,E[p^\infty])$
that sits in an exact sequence
$$
0 \longrightarrow E(F) \otimes \Qp/\Zp \longrightarrow \Sel_p(E_{/F})
     \longrightarrow \Sh(E_{/F})[p^\infty] \longrightarrow 0
$$
where $\Sh(E_{/F})$ is the Shafarevich-Tate group of $E$ over $F$.
Also write
$$
\Selp{F} = \Hom(\Sel_p(E_{/F}),\Qp/\Zp)\otimes\Qp
$$
for the tensor product of $\Qp$ with the Pontrjagin dual
of the Selmer group.

The following theorem is proved using techniques which go back to
\cite{mazuritav}; see \cite{greenberg} and \cite{perrinrioubsmf}
Lemme 5, \S2.2.

\begin{thm}[Control Theorem]
\label{control}
Suppose $K \subset F \subset \Ktwo$.
\begin{enumerate}
\item
The natural restriction map
$H^1(F,E[p^\infty]) \to H^1(\Ktwo,E[p^\infty])$ induces an isomorphism
$
\Selp{\Ktwo} \otimes_{\bLambda} \Lambda_F \isom \Selp{F}.
$
\item
There is a canonical isomorphism
$U(F) \isom \Hom_{\Lambda_F}(\Selp{F},\Lambda_F)$.
\end{enumerate}
\end{thm}

\begin{conj}[Two-variable main conjecture \cite{mazuricm, perrinrioubsmf}]
\label{twomc}
The two-variable $p$-adic $L$-function $\L$ generates the ideal
$\cp_\bLambda(\Selp{\Ktwo})$ of $\bLambda$.
\end{conj}

Restricting the two-variable main conjecture to the cyclotomic
and anticyclotomic lines leads to the following
``one-variable'' conjectures originally formulated in \cite{mazursd}
and \cite{perrinrioubsmf}, respectively.
Let $L'$ denote the image of $\L$ in
$\Gcycl\otimes_{\Zp}\Lambda_{\anti}$ as in Conjecture \ref{lambdagz},
and $\Selp{\Kanti}_{\tors}$ the $\Lambda_{\anti}$-torsion
submodule of $\Selp{\Kanti}$).

\begin{conj}[Cyclotomic and anticyclotomic main conjectures]~
\label{onemc}
\begin{enumerate}
\item
$L_{\cycl}$ generates the ideal $\cp_{\Lambda_{\cycl}}(\Selp{\Kcycl})$
of $\Lambda_{\cycl}$.
\item
$L'$ generates
$\Gcycl \otimes \cp_{\Lambda_{\anti}}(\Selp{\Kanti}_{\tors})$
inside $\Gcycl \otimes \Lambda_{\anti}$.
\end{enumerate}
\end{conj}

\begin{rem}
Using Euler systems, Kato \cite{kato} and Howard \cite{howard}
have proved (under some mild additional hypotheses)
divisibilities related to the
cyclotomic and anticyclotomic main conjectures, respectively, namely
$$
L_{\cycl} \in \cp_{\Lambda_{\cycl}}(\Selp{\Kcycl}), \quad
\cp_{\Lambda_{\anti}}(\U/\H)^2 \subset
    \cp_{\Lambda_{\anti}}(\Selp{\Kcycl}_{\tors})
$$
(note that Conjectures \ref{heegnerlfunction} and \ref{lambdagz}
predict that $\Gcycl \otimes \cp_{\Lambda_{\anti}}(\U/\H)^2 = L'\Lambda_{\anti}$).
\end{rem}

\begin{conj}[Two-variable $p$-adic BSD conjecture]
\label{twoBSD}
Let $r = \rk{E(K)}$. The two-variable
$p$-adic $L$-function $\L\in\bLambda$ is contained $\I^r$ and
$$
\L \;\equiv\; c(\trivchar) \#(\Sh(E_{/K})) \prod_v c_v \cdot R_p(E,K) \pmod{\I^{r+1}}
$$
where $c(\trivchar)$ is the rational
number in the interpolation formula \eqref{const} for the trivial character,
$\Sh(E_{/K})$ is the Shafarevich-Tate group of $E$ over $K$,
and the $c_v$ are the Tamagawa factors in
the (usual) Birch and Swinnerton-Dyer conjecture for $E$ over $K$.
\end{conj}

\section{Orthogonal $\bLambda$-modules} \label{orthogonal}\setzero
\vskip-5mm \hspace{5mm }

In this final section we introduce a purely algebraic template
which, when it ``fits'', gives rise to many of the properties conjectured
in the previous sections.

Keep the notation of the previous sections.  In particular
$\tau : \bLambda \to \bLambda$ is the involution of $\bLambda$ induced by
complex conjugation on $K$, and if
$V$ is a $\bLambda$-module, then $V^{\opp}$ denotes $V$ with
$\bLambda$-module structure obtained by composition with $\tau$.
Let $V^* = \Hom_\bLambda(V,\bLambda)$.
If $V$ is a free $\bLambda$-module of rank $r$, then $\det_\bLambda(V^{\tau})$
will denote the $r$-th exterior power of $V$ and a
{\bf $\tau$-gauge} on $V$ is a $\bLambda$-isomorphism between the free
$\bLambda$-modules of rank one
$$
\textstyle
t_V \;:\; \det_\bLambda(V^*) \cong \det_\bLambda(V^{\opp})
$$
or equivalently an isomorphism
$\det_\bLambda(V) \otimes \det_\bLambda(V^{\opp}) \isom \bLambda$.

By an {\bf orthogonal} $\bLambda$-module 
we mean a free $\bLambda$-module $V$ with semi-linear
$\tau$-structure endowed with a $\tau$-gauge $t_V$ and
a $\bLambda$-bilinear $\tau$-Hermitian pairing (Definition \ref{herm})
$$
\pi: V \otimes_\bLambda V^{\opp} \longrightarrow \bLambda.
$$
Viewing $\pi$ as a $\bLambda$-linear map $V^{\opp} \to V^*$, the composition
$$
\textstyle
t_V \circ \det_\bLambda(\pi) \;:\; \det_\bLambda(V^{\opp})
    \longrightarrow \det_\bLambda(V^*) \longrightarrow \det_\bLambda(V^{\opp})
$$
must be multiplication by an element $\disc(V) \in \bLambda$ that we call
the {\bf discriminant} of the orthogonal $\bLambda$-module $V$.
We further assume that $\disc(V) \ne 0$, and we define $M = M(V,\pi)$ to be
the cokernel of the (injective) map $\pi : V^{\opp} \to V^*$, so we have
\begin{equation}
\label{Mdef}
0 \longrightarrow V^{\opp} \longrightarrow V^* \longrightarrow M
    \longrightarrow 0.
\end{equation}

If $K \subset F \subset \Ktwo$, recall that
$\I_F = \ker \{\bLambda \twoheadrightarrow \Lambda_F\}$ and define
$$
V(F) := \{x \in V : \pi(x,V^{\opp}) \subset \I_F\}/\I_F V
    = \ker\{V \otimes_\bLambda \Lambda_F
    \map{\pi\otimes 1} (V^{\opp})^* \otimes_\bLambda \Lambda_F\}
$$
and similarly
$V^{\opp}(F) := \ker\{V^{\opp} \otimes \Lambda_F \to V^* \otimes \Lambda_F\}$.
Any lift $\tilde\tau$ of $\tau$ to $\Gal(\Ktwo/\Q)$ induces an isomorphism
$V(F) \to V^{\opp}(F)$.
From \eqref{Mdef} we obtain
\begin{equation}
\label{pre}
0 \longrightarrow V^{\opp}(F) \longrightarrow V^{\opp} \otimes_{\bLambda} \Lambda_F
    \longrightarrow V^* \otimes_{\bLambda} \Lambda_F
    \longrightarrow M \otimes_{\bLambda} \Lambda_F \longrightarrow 0
\end{equation}
and (applying $\Hom(\;\cdot\;,\Lambda_F)$ and using the Hermitian property of $\pi$)
\begin{equation}
\label{MmodI}
V(F) \cong \Hom_{\Lambda_F}(M \otimes_{\bLambda} \Lambda_F, \Lambda_F).
\end{equation}

We have an induced pairing
$$
\pi_F : V^{\opp}(F) \otimes_{\Lambda_F} V(F) \longrightarrow \I_F/\I_F^2,
$$
which we call the $F$-derived pairing.  If $F$ is stable under complex
conjugation then $V^{\opp}(F)$ is canonically isomorphic to
$V(F)^{\opp}$ and $\pi_F$ is $\tau$-Hermitian.

Now suppose $F = \Kanti$.  By \eqref{MmodI}, $V(\Kanti)$ is free over $\Lambda_{\anti}$.
Applying the determinant functor to \eqref{pre}, the $\tau$-gauge $t_V$ induces an
isomorphism
$$
\textstyle
\det_{\Lambda_{\anti}} V(\Kanti)^{\opp} = \det_{\Lambda_{\anti}} V^{\opp}(\Kanti)
    \isom \Hom(\det_{\Lambda_{\anti}}(M \otimes_{\bLambda} \Lambda_{\anti}), \Lambda_{\anti}).
$$
If $V(\Kanti)$ has rank one over $\Lambda_{\anti}$, then $V(\Kanti)$ contains a
unique maximal $\tau$-stable submodule $H$ such that the map
\begin{multline*}
\textstyle
V(\Kanti)^{\opp} \isom
    \Hom(\det_{\Lambda_{\anti}}(M \otimes \Lambda_{\anti}), \Lambda_{\anti}) \\
    \supset \Hom(M \otimes \Lambda_{\anti}, \Lambda_{\anti}) \cong V(\Kanti)
\end{multline*}
sends $H^{\opp}$ into $H$.  (Namely, $H = \cA V(\Kanti)$ where
$\cA$ is the largest ideal of $\Lambda_{\anti}$ such that $\cA^\tau = \cA$ and
$\cA^2 \subset \cp_{\Lambda_{\anti}}(M \otimes \Lambda_{\anti})_{\tors}$.)

Recall that $\Sel_p(E_{/F})$ denotes the $p$-power Selmer group of $E$ over $F$ and
$\Selp{F} = \Hom(\Sel_p(E_{/F}),\Qp/\Zp)\otimes\Qp$.

\begin{prop}
\label{inducedmaps}
With notation as above, suppose that $V$ is an orthogonal $\bLambda$-module
and $\varphi_V : M \isom \Selp{\Ktwo}$ is an isomorphism.  Then for every
extension $F$ of $K$ in $\Ktwo$, $\varphi_V$ induces an isomorphism
$$
V(F) \isom U(F)
$$
where $U(F)$ is the universal norm module defined in \S\ref{universalnorms}
\end{prop}

\noindent{\bf proof.} This follows directly from Theorem
\ref{control} and \eqref{MmodI}.

\begin{defn}
\label{orgdef}
\renewcommand{\labelenumi}{\theenumi}
\renewcommand{\theenumi}{(\alph{enumi})}
\it
We say that the orthogonal $\bLambda$-module $V$
{\bf organizes the anticyclotomic arithmetic of $(E,K,p)$} if
the following three properties hold.
\begin{enumerate}
\item
$\disc(V) = \L$, the two-variable $p$-adic $L$-function of $E$.
\item
There is an isomorphism
$
\varphi_V : M \isom \Selp{\Ktwo}.
$
\item
The isomorphism
$V(\Kanti) \cong \U$ of Proposition \ref{inducedmaps} identifies
$H \subset V(\Kanti)$ with the Heegner submodule $\H \subset \U$, and identifies the
$\Kanti$-derived pairing with the canonical $p$-adic height pairing into
$\I_{\Kanti}/\I_{\Kanti}^2 \cong \Gcycl \otimes \Lambda_{\anti}$.
\end{enumerate}\rm
\end{defn}

\begin{question}\rm
Given $E$, $K$, and $p$ satisfying our running hypotheses, is there
an orthogonal $\bLambda$-module $V$ that
organizes the anticyclotomic arithmetic of $(E,K,p)$?
\end{question}

If one is not quite so (resp., much more) optimistic one could formulate
an analogous question with the ring $\bLambda$
replaced by the localization of $\bLambda$ at $\I$ (resp., with $\bLambda$ replaced
by $\Zp[[\bGamma]]$).

\begin{question}\rm
If $V$ is an orthogonal $\bLambda$-module $V$ which
organizes the anticyclotomic arithmetic of $(E,K,p)$, then for every finite extension
$F$ of $K$ in $\Ktwo$, we have an isomorphism $E(F) \otimes \Qp = U(F) \cong V(F)$
as in Proposition \ref{inducedmaps}, a $p$-adic height pairing on $E(F) \otimes \Qp$, and
the $F$-derived pairing on $V(F)$.  How are these pairings related?
\end{question}

When $F = \Kanti$
condition (c) says that the two pairings are the same, but it seems that in general
they cannot be the same for finite extensions $F/K$.

\begin{thm}
\label{organized}
Suppose that there is an orthogonal $\bLambda$-module $V$ that
organizes the anticyclotomic arithmetic of $(E,K,p)$.  Then
Conjectures \ref{twomc} (the 2-variable main conjecture), and
\ref{onemc}(i) (the cyclotomic main conjecture)
hold.

If further the induced pairing
$V(\Kanti) \otimes V(\Kanti)^{\opp} \to \Gcycl \otimes \Lambda_{\anti}$
is surjective, then Conjectures \ref{heightconj} (the Height Conjecture),
\ref{heegnerlfunction},
\ref{lambdagz} (the $\Lambda$-adic Gross-Zagier conjecture), and
\ref{onemc}(ii) (the anticyclotomic main conjecture) also hold.
\end{thm}

\begin{proof}[Brief outline of the proof of Theorem
\ref{organized}] Since $\disc(V)$ is a generator of
$\cp_{\bLambda}(M)$, the two-variable main conjecture follows
immediately from (a) and (b) of Definition \ref{orgdef}. The
cyclotomic main conjecture follows from the two-variable main
conjecture.

Now suppose that the induced pairing
$V(\Kanti) \otimes V(\Kanti)^{\opp} \to \Gcycl \otimes \Lambda_{\anti}$
is surjective.  By (c) of Definition \ref{orgdef} this is equivalent
to the Height Conjecture, which in turn is equivalent to
Conjecture \ref{heegnerlfunction}.

Howard proved in \cite{howard} that $\Selp{\Kanti}$ is pseudo-isomorphic to
$\Lambda_{\anti} \oplus B^2$ where $B$ is a $\tau$-stable torsion $\Lambda_{\anti}$-module.
By Theorem \ref{control}(i) the same is true of $M \otimes \Lambda_{\anti}$, and
so the remark at the end of the definition of $H$ shows that $H = \cp(B)V(\Kanti)$.
Using \eqref{pre}, \eqref{MmodI}, and our assumption that the induced pairing is
surjective, one can show that
the image of $\L$ in $\I_{\Kanti}/\I_{\Kanti}^2$ generates
$\cp(B)^2\I_{\Kanti}/\I_{\Kanti}^2$.
The $\Lambda$-adic Gross-Zagier conjecture and the anticyclotomic main conjecture
follow from these facts and (c).
\end{proof}

\label{lastpage}

\end{document}